\documentclass{amsart}
\usepackage{upref}
\usepackage{graphicx}
\usepackage{pstricks}
\usepackage{pst-plot}

\newtheorem{theorem}{Theorem}[section]

\newtheorem{pro}[theorem]{Proposition}

\theoremstyle{definition}

\theoremstyle{remark}

\numberwithin{equation}{section}

\newcommand{\vphi}{\varphi}
\newcommand{\set}[1]{\left\{#1\right\}}

\newcommand{\CC}{{\mathbb C}}

\newcommand{\DD}{{\mathbb D}}

\newcommand{\ca}[1]{{\mathcal #1}}

\begin{document}

\title{Composition Operators on the Dirichlet space and related problems}
\author{Gerardo A. Chac\'{o}n}
\address{Universidad de los Andes \\N\'{u}cleo T\'{a}chira}
\email{gchacon@cantv.net, gchacon@ula.ve}
\thanks{}

\author{Gerardo R. Chac\'{o}n}
\address{Universidad de los Andes\\Facultad de Humanidades y Educaci\'{o}n}
\email{grchacon@ula.ve}
\thanks{}

\author{Jos\'{e} Gim\'{e}nez}
\address{Universidad de los Andes\\Facultad de Ciencias}
\email{jgimenez@ula.ve}
\thanks{The authors are partially supported by a grant of CDCHT-ULA, Venezuela.}

\subjclass[2000]{Primary: 47B33; Secondary 47B38, 47A16}
\date{\today}

\begin{abstract}
In this paper we investigate the following problem: when a bounded analytic function $\vphi$ on the unit
disk $\DD$, fixing $0$, is such that $\set{\vphi^n:n=0,1,2,\dots}$ is orthogonal in $\ca{D}$?, and
consider the problem of characterizing the univalent, full self-maps of $\DD$ in terms of the norm of
the composition operator induced. The first problem is analogous to a celebrated question asked by W.
Rudin on the Hardy space setting that was answered recently (\cite{Bi} and \cite{Su}). The second
problem is analogous to a problem investigated by J. Shapiro in \cite{Sh} about characterization of
inner functions in the setting of $H^2$.
\end{abstract}

\maketitle

Let $\DD$ denote the unit disk in the complex plane. By a {\em self- map} of $\DD$ we mean an analytic
map such that $\vphi(\DD)\subset \DD$. The {\em composition operator} induced by $\vphi$ is the linear
transformation $C_\vphi$ defined by $C_\vphi(f)=f\circ \varphi$ in the space of the holomorphic
functions on $\DD$.

The composition operators have been studied in many settings, and in particular in functional Banach
spaces (cf. the books \cite{CM}, \cite{Sh1}, the survey of recent developments \cite{Ja}, and the
references therein). Recall that a functional Banach space  is a Banach space of analytic functions (on
the disk or other domains of $\CC$ or $\CC^n$) where the evaluation functionals are continuous. The goal
of this theory is to obtain characterizations of operator-theoretic properties of $C_\vphi$ by
function-theoretic properties of the symbol $\vphi$. Conversely, operator-theoretic properties of
$C_\vphi$ could  suggest, or help to understand certain phenomena about function-theoretic properties of
$\vphi$.

Particular instances of functional Banach spaces are the Hardy space
$H^2$, and  the Bergman space $A^2$ of the unit disk. In these
spaces, as a consequence of Littlewood's Subordination Principle,
every self-map of $\DD$ induces a bounded composition operator. A
very interesting setting for studying such operators is the
Dirichlet space. Recall that if
$dA(z)=\frac{1}{\pi}dx\,dy=\frac{1}{\pi}r\,dr\,d\theta$,
($z=x+iy=re^{i\theta}$) denotes the normalized area Lebesgue measure
on $\DD$, the {\em Dirichlet space} $\ca{D}$ is the Hilbert space of
analytic functions in $\DD$ with a square integrable derivative,
with the norm given by
$$
\|f\|_{\ca{D}}^2=|f(0)|^2+\int_{\DD}|f'(z)|^2\,dA(z).
$$

It is well known that $\ca{D}$ is a functional Hilbert space, and
for each $w\in \DD$ the function
$$
K_w(z)=1+\log\frac{1}{1-\overline{w}z},
$$
is  the reproducing kernel at $w$ in the Dirichlet space, that is,
for $f\in \ca{D}$ we have $\langle f,K_w\rangle_{\ca{D}}=f(w)$. It
is easy to see that $\|K_w\|_{\ca{D}}^2=\log\frac{1}{1-|w|^2}$.

A self-map of $\DD$ does not induce, necessarily a bounded
composition operator on $\ca{D}$. An obvious necessary condition for
it is that $\vphi=C_\vphi z\in \ca{D}$ which is not always the case.
Actually this condition is not sufficient. A necessary and
sufficient condition in order to $\vphi$ to induce a bounded
composition operator on $\ca{D}$ is given in terms of counting
functions and Carleson measures (see \cite{JM} and the references in
this paper).

Recall that the {\em counting function} $n_\vphi(w)$, $w\in \DD$,
associated to $\vphi$ is defined as the cardinality of the set
$\set{z\in \DD:\vphi(z)=w}$ when the latter is finite and understood
as the symbol $\infty$ otherwise, with the usual rules of
arithmetics holding  in relation to the Lebesgue integral.

We will make use of a change of variable formula for non-univalent
functions: Suppose $\vphi:\DD\to\DD$ is a non-constant analytic
function with counting function $n_\vphi(w)$, if
$f:\DD\to[0,\infty)$ is any Borel function, then
$$
\int_{\DD} f(\vphi(z))|\vphi'(z)|^2\,dA(z)=\int_{\DD}
f(w)n_\vphi(w)\,dA(w).
$$
This formula is a particular instance of the general change of
variable formula in \cite[Th. 2.32]{CM} (see also \cite{Fe}). In
particular we obtain, as noted in \cite{Fe},  that
$\int_{\DD}|\vphi'(z)|^2\,dA(z)=\int_{\DD}n_\vphi(w)\,dA(w)$. So,
$\vphi$ is in the Dirichlet space if and only if its counting
function is an $L^1$ function.

In two recent papers, \cite{MV1} and \cite{MV2}, M. Mart\'{\i}n and D. Vukoti\'{c},  studied composition
operators on the Dirichlet space. In this note, based on results in those works, we consider related
questions. In Section 1, we investigate the analogous on Dirichlet spaces to a problem proposed by W.
Rudin in the context of Hardy spaces: When a bounded analytic functions $\vphi$  on the unit disk $\DD$
fixing $0$ is such that $\set{\vphi^n:n=0,1,2,\dots}$ is orthogonal in $\ca{D}$?, and in  Section 2 we
consider the problem of characterizing the univalent, full self-maps of $\DD$ in terms of the norm of
the composition operator induced. This problem, is analogous to the question asked and answered by J.
Shapiro in \cite{Sh} about inner functions in the  $H^2$ setting.

We write $\ca{D}_0$ for the subspace of $\ca{D}$ of those function in $\ca{D}$ vanishing in $0$, and use
the notation $\|C_\vphi:\ca{H}\to\ca{H}\|$ in order to denote the norm of the composition operator
induced on the space $\ca{H}$.

\section{Orthogonal functions in the Dirichlet space.}
The problem of describing  the isometric composition operators
acting in Hilbert spaces of analytic functions has been studied in
several settings. Namely, it was proved by Nordgren in \cite{No}
that the composition operator $C_\vphi$ induced on $H^2$ by $\vphi$,
a holomorphic self-map  of the unit disk, is an isometry on $H^2$ if
and only if $\vphi(0)=0$ and $\vphi$ is an inner function (see also
\cite[p. 321]{CM}. In $A^2$ it is an straightforward consequence of
the Schwarz Lemma that $\vphi$ induces an isometric composition
operator if and only if it is a rotation.

Recently, M. Martin and D. Vukoti\'{c} showed in \cite{MV2} that in
$\ca{D}$, the Dirichlet space in the unit disk, the isometric
composition operators are those induced by univalent full maps of
the disk into itself that fixes the origin. Recall that a self-map
of $\DD$ is said a {\em full map} if $A[\DD\setminus\vphi(\DD)]=0$.

W. Rudin in 1988 (at an MSRI conference) proposed  the following
problem: If $\vphi$ is a bounded analytic on the unit disk $\DD$
such that $\set{\vphi^n:n=0,1,2,\dots}$ is orthogonal in $H^2$, does
$\vphi$ must be a constant multiple of an inner function? C.
Sunbberg \cite{Su} and C. Bishop \cite{Bi} solved independently the
problem. In fact, they show that there exists a function $\vphi$
such that $\vphi$ is not an inner function and $\set{\vphi^n}$ is
orthogonal in $H^2$.

As  asserted by M. Mart\'{\i}n y D. Vukoti\'{c} in \cite{MV2}, their
characterization of the isometric composition operators acting on
$\ca{D}$ can be interpreted as  follows: the univalent full maps of
the disk that fix the origin are the Dirichlet space counterpart of
the inner functions that fix the origin for the composition
operators on $H^2$. We  propose the following question: When a
bounded analytic function $\vphi$  on the unit disk $\DD$ fixing $0$
is such that $\set{\vphi^n:n=0,1,2,\dots}$ is orthogonal in
$\ca{D}$? Recall that a bounded analytic function on $\DD$ is not
necessarily in $\ca{D}$, then we assume in this context that $\vphi$
is in $\ca{D}$ (and therefore, since $\ca{D}\cap H^\infty$ is an
algebra, that $\set{\vphi^n}$ is in $\ca{D}$).

We are going to answer this question in the case when $n_\vphi$ is essentially bounded, that is, there
is a constant $C$ so that $n_\vphi(w)\leq C$ for all $w$ except those in a set of area zero. Our result
is analogous to a characterization given by P. Bourdon in \cite{Bo} in the context of $H^2$: the
functions that satisfy the hypotheses of the Rudin's problem are characterized as those maps $\vphi$
such that their Nevanlinna counting function $N_\vphi$ is essentially radial. Our assumption that
$n_\vphi$ is essentially bounded is clearly stronger that assuming that $\vphi$ is only in the Dirichlet
space and it possibly can be relaxed. The proof relies in the techniques of the proof given in \cite
{Bo}.

\begin{theorem}
Let $\vphi$ be a self-map on $\DD$ fixing $0$. The set
$\set{\vphi^n:n=0,1,2,\dots}$ is orthogonal in $\ca{D}$ if and only
if there is a function $g:[0,1)\to [0,\infty)$ such that for almost
every $r\in [0,1)$, $n_\vphi(re^{i\theta})=g(r)$ for almost every
$\theta\in [0,2\pi]$ (this is, $n_\vphi$ is essentially radial).
\end{theorem}
\begin{proof}
Suppose that $n_\vphi$ is essentially radial. Let $n>m$ be
nonnegative integers. We have
\begin{align*}
\langle \vphi^n,\vphi^m\rangle_{\ca{D}}
    &=nm\int_{\DD}\vphi(z)^{n-1}\overline{\vphi(z)^{m-1}}|\vphi'(z)|^2\,dA(z)\\
    &=nm\int_{\DD}w^{n-1}\overline{w^{m-1}}n_\vphi(w)\,dA(w)\\
    &=nm\int_0^1
    r^{n+m-1}\left[\frac{1}{\pi}\int_0^{2\pi}e^{i(n-m)\theta}n_\vphi(re^{i\theta})\,d\theta\right]\,dr\\
    &=nm\int_0^1
    r^{n+m-1}g(r)\left[\frac{1}{\pi}\int_0^{2\pi}e^{i(n-m)\theta}\,d\theta\right]\,dr\\
    &=0.
\end{align*}
Conversely, if $\set{\vphi^n:n=0,1,2\dots}$ is orthogonal in
$\ca{D}$. Let $k$ be an arbitrary positive integer. For each integer
$n>k$, we have
\begin{align*}
0=\langle \vphi^n,\vphi^{n-k}\rangle_{\ca{D}}
    &=n(n-k)\int_{\DD}\vphi(z)^{n-1}\overline{\vphi(z)^{n-k-1}}|\vphi'(z)|^2\,dA(z)\\
    &=n(n-k)\int_{\DD}w^{n-1}\overline{w^{n-k-1}}n_\vphi(w)\,dA(w)\\
    &=n(n-k)\int_0^1
    r^{2n-k-1}\left[\frac{1}{\pi}\int_0^{2\pi}e^{ik\theta}n_\vphi(re^{i\theta})\,d\theta\right]\,dr.
\end{align*}
The functions
$f_k(r):=\int_0^{2\pi}e^{ik\theta}n_\vphi(re^{i\theta})\,d\theta$
are in $L^2[0,1]$  since $n_\vphi$ is essentially bounded (it is the
only instance of this hypothesis) and the precedent equation says
that they are orthogonal in $L^2[0,1]$ to $\set{r\mapsto
r^{2n-k-1}:n>k}$. By an slight variation of M\"{u}ntz-Szasz Theorem (cf.
\cite{Bo}), the linear span of this set is dense in $L^2[0,1]$, and
so $f_k(r)=0$ for almost every $r\in [0,1]$. Taking complex
conjugates, we see that
$\int_0^{2\pi}e^{ij\theta}n_\vphi(re^{i\theta})\,d\theta=0$ for all
$j\neq 0$, and almost every $r\in [0,1]$. Thus that $\theta\mapsto
n_\vphi(re^{i\theta})$ is essentially constant for almost every $r$.
\end{proof}

The following Proposition describes the self-maps of $\DD$ that
share the properties in the condition of the the previous Theorem.

\begin{pro}\label{cor}
Suppose that $\vphi$ is a self-map with counting function
essentially bounded, and essentially radial. Then $\vphi$ is a
constant multiple of a full self-map of $\DD$.
\end{pro}
\begin{proof}
Suppose that $\vphi$ is not constant. If the range  of $\vphi$
contains a point in the circle $S_r=\set{re^{i\theta}:\theta\in
[0,2\pi]}$, $\vphi(\DD)$ contains an arc because this is an open
subset of $\DD$. In this arc $n_\vphi\geq 1$,  and so the range of
$\vphi$ may omit only a $\theta$-zero-measure subset of $S_r$
because $n_\vphi$ is essentially constant on $S_r$.

Thus the range of $\vphi$ contain almost every point in the disk
$\set{z:|z|<\|\vphi\|_\infty}$.
\end{proof}

\section{What do composition operators know about full mappings?}
In the Hardy space, J. Shapiro \cite{Sh} has characterized, in terms of their norms, those composition
operators $C_\vphi$ whose symbol is an {\em inner function.} In fact, J. Shapiro showed:
\begin{enumerate}
    \item If $\vphi(0)=0$ then $\vphi$ is inner if and only if
    $\|C_\vphi:H_0^2\to H_0^2\|=1$, where $H_0^2$ is the subspace of
    functions in $H^2$ what vanish at $0$, and
    \item If $\vphi(0)\neq 0$ then $\vphi$ is inner if and only if
    $\|C_\vphi:H^2\to H^2\|=\sqrt{\frac{1+|\vphi(0)}{1-|\vphi(0)|}}$.
\end{enumerate}
We are going to investigate the analogous questions on the Dirichlet space.

In \cite{MV1} M. Mart\'{\i}n and D. Vukoti\'{c} calculate the norm of the composition operator $C_\vphi$
induced on $\ca{D}$ by a univalent full map $\vphi$ of $\DD$. They obtain
\begin{equation}\label{Norm}
\|C_\vphi:{\ca D}\to{\ca D}\|=\sqrt{\frac{L+2+\sqrt{L(4+L)}}{2}},
\end{equation}
where $L=\log\frac{1}{1-|\vphi(0)|^2}$, and show that it is an upper bound on the norms of composition
operators acting on the Dirichlet space induced by univalent symbols.

The results in \cite{Sh} and the assertion in \cite{MV2}, mentioned previously, that the univalent full
maps of the disk that fix $0$ are the Dirichlet space counterpart of the inner functions that fix the
origin for the composition operators on $H^2$, lead us to investigate if the equality in the equation
(\ref{Norm}) characterizes the univalent full maps of the disk inside the univalent self-maps of $\DD$.

In addition, the main result in \cite{MV2} says that $\vphi(0)=0$ and $\vphi$ is a univalent full
self-map of the disk if and only if $C_\vphi$ is an isometry on $\ca{D}$, and hence on $\ca{D}_0$, so in
particular its restriction to $\ca{D}_0$ has norm $1$. Is the converse true?

It is easy to see that this is not true. In fact, let $\vphi_t$,
$t\geq 1$, be the linear fractional transformation given by
$$
\varphi_t(z)=\frac{2z}{(1-t)z+(1-z)},\quad z\in \DD.
$$
We easily see  that $\varphi_t(\DD)\subset\DD$, $\vphi_t(0)=0$,
$\vphi_t(1)=1$, and $\vphi_t(-1)=-1/t$ (see figure). If $t>1$
clearly $\vphi_t$ is not full, but a calculation in \cite[Cor.
6.1]{02} shows that $\|C_\vphi:\ca{D}_0\to \ca{D}_0\|=1$ when
$\vphi$ is a linear fractional self-map of $\DD$ with a boundary
fixed point.

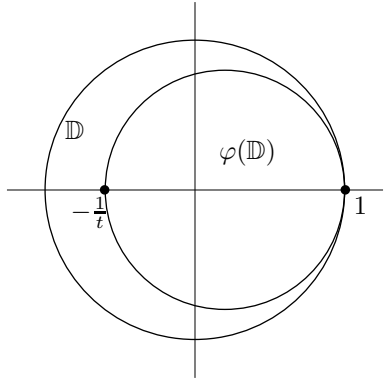
\begin{center}
\begin{pspicture}(5,5)
\psaxes[linewidth=0.3pt,labels=none,
ticks=none]{-}(2.5,2.5)(0,0)(5,5)
\pscircle[linewidth=.5pt](2.5,2.5){2}
\pscircle[linewidth=.5pt](2.9,2.5){1.6}
\rput(3.2,3){\shortstack{$\vphi(\DD)$}}
\rput(.9,3.3){\shortstack{$\DD$}} \rput(4.7,2.3){\shortstack{$1$}}
\rput(1.1,2.2){\shortstack{$-\frac{1}{t}$}}
\psdots(4.5,2.5)(1.3,2.5)
\end{pspicture}
\end{center}

Nevertheless, we have the following results, analogous to the results in \cite{Sh}.
\begin{theorem}\label{T1}
Suppose $\varphi$ is a univalent, holomorphic self-map of $\DD$,
with $n_\vphi$ essentially radial and $\vphi(0)=0$. Then $\vphi$ is
a full map if and only if
$$\|C_\vphi:\ca{D_0}\to\ca{D_0}\|=1.$$
\end{theorem}
\begin{proof}
We saw before one direction. For the converse, suppose that $\vphi$
is  a univalent holomorphic self-map of $\DD$, with $n_\vphi$
essentially radial, $\vphi(0)=0$, and that $\vphi$ is not a full
map.

We are going to show that the restriction of $C_\vphi$ to $\ca{D}_0$
has norm $<1$. We have that $\vphi(\DD)$ is contained in the disk
$D(0,\rho)=\set{z:|z|<\|\vphi\|_\infty=\rho}$ with
$A[\vphi(\DD)\setminus D(0,\rho)]=0$ and $0<\rho<1$ (cf. proof of
Proposition \ref{cor}.)

We write
$$
g(r):=\frac{1}{\pi}\int_0^{2\pi} |f'(re^{i\theta})|^2\,d\theta,
$$
and since $|f'|^2$ is subharmonic in $\DD$ then $g$ is monotone
increasing for $0\leq r<1$. The change of variable formula gives
\begin{align*}
\|C_\vphi\|_{\ca{D}}^2&=\int_{\DD}|f'(\vphi(z))|^2\,|\vphi'(z)|^2\,dA(z)\\
    &=\int_{\vphi(\DD)}|f'(w)|^2\,dA(w)\\
    &=\int_0^\rho g(r)\,r\,dr.
\end{align*}

and so:
\begin{align*}
\|f\|_{\ca{D}}^2=\int_{\DD}|f'(w)|^2\,dA(w)&=\int_0^\rho
g(r)\,r\,dr+\int_\rho^1 g(r)\,r\,dr\\
    &\geq \int_0^\rho g\,dr+\frac{1-\rho^2}{2}g(\rho)\\
    &=\int_0^\rho g\,dr+\frac{(1-\rho^2)/2}{\rho^2/2}(\rho^2/2)g(\rho)\\
    &\geq \int_0^\rho g\,dr+\frac{(1-\rho^2)/2}{\rho^2/2}\int_0^\rho g(r)\,r\,dr\\
    &=\left(1+\frac{(1-\rho^2)/2}{\rho^2/2}\right)\int_0^\rho
    g(r)\,r\,dr\\
    &=\left(1+\frac{(1-\rho^2)/2}{\rho^2/2}\right)\|C_\vphi\|_{\ca{D}}^2,
\end{align*}
for each $f\in \ca{D}_0$. It yields the desired result: the
restriction of $C_\vphi$ to $\ca{D}_0$ has norm $\leq
\nu=\left(1+\frac{(1-\rho^2)/2}{\rho^2/2}\right)^{-1/2}<1$.
\end{proof}

In the next theorem, we consider the case $\vphi(0)\neq 0$. The proof follows nearly the one in
\cite[Th. 5.2]{Sh}).

\begin{theorem}
Suppose $\varphi$ is a univalent, holomorphic self-map of $\DD$ with
$n_\vphi$ essentially radial and $\vphi(0)\neq 0$. Then $\vphi$ is a
full map if and only if
$$
\|C_\vphi:\ca{D}\to\ca{D}\|=\sqrt{\frac{L+2+\sqrt{L(4+L)}}{2}},
$$
where $L=\log 1/(1-|\vphi(0)|^2)$.
\end{theorem}
\begin{proof}
The necessity is part of \cite[Th. 1]{MV1}. For the  converse,
suppose that $\vphi$ is a univalent, holomorphic self-map of $\DD$
with $n_\vphi$ essentially radial, such that $\vphi(0)=p\neq 0$, and
$\vphi$ is not a full map. We want to show that the norm of
$C_\vphi$ is strictly less that
$\sqrt{\frac{L+2+\sqrt{L(4+L)}}{2}}$, where $L=\log 1/(1-p^2)$.

For this  we consider $\alpha_p$, the standard automorphism of $\DD$
that interchanges $p$ with the origin, this is
$$
\alpha_p:=\frac{p-z}{1-\overline{p}z},\quad z\in \DD.
$$
We write $\vphi_p=\alpha_p\circ \vphi$, which is $0$ in the origin.
Since this function is a univalent, self map of $\DD$ with counting
function essentially radial, but it is not full, the Theorem
\ref{T1} affirms  that the restriction of the operator $C_{\vphi_p}$
to $\ca{D}_0$ has norm $\nu<1$.

Because  $\alpha_p$ is self-inverse, $\vphi=\alpha_p\circ\vphi_p$,
and so, for each $f\in \ca{D}$:
$$
C_\varphi f=C_{\vphi_p}(f\circ\alpha_p)=C_{\vphi_p}f+f(p),
$$
where $g=f\circ \alpha_p-f(p)$.

The function $C_{\vphi_p} g$ belong to $\ca{D}_0$ and thus:
\begin{align}\label{Eq2}
\|C_\phi f\|_{\ca{D}}=&\|C_{\vphi_p}g\|_{\ca{D}}^2+|f(p)|^2\\ \notag
    &\leq \nu^2\|g\|_{\ca{D}}^2+|f(p)|^2\\ \notag
    &=\nu^2\|(C_{\alpha_p}f)-f(p)\|_{\ca{D}}^2+|f(p)|^2.
\end{align}
Since $\langle h,1\rangle_{\ca{D}}=h(0)$ for each $h\in \ca{D}$,
$$
\langle
C_{\alpha_p}f,f(p)\rangle_{\ca{D}}=\overline{f(p)}C_{\alpha_p}f(0)=|f(p)|^2,
$$
and we obtain
\begin{align*}
\|(C_{\alpha_p}f)-f(p)\|_{\ca{D}}^2&=\|C_{\alpha_p}f\|_{\ca{D}}^2-2\Re
\langle C_{\alpha_p}f,f(p)\rangle_{\ca{D}}+|f(p)|^2\\
    &=\|C_{\alpha_p}f\|_{\ca{D}}^2-2|f(p)|^2+|f(p)|^2\\
    &=\|C_{\alpha_p}f\|_{\ca{D}}^2-|f(p)|^2.
\end{align*}
This identity and the Equation (\ref{Eq2}) yield,
$$
\|C_{\alpha}f\|_{\ca{D}}^2\leq
\nu^2\|C_{\alpha_p}f\|_{\ca{D}}^2+(1-\nu^2)|f(p)|^2.
$$
We know from \cite[Th. 1]{MV2} that
$\|C_{\alpha_p}:\ca{D}\to\ca{D}\|=(L+2+\sqrt{L(4+L)})/2$, and we
have the following estimate for $|f(p)|$:
$$
|f(p)|\leq \|f\|_{\ca{D}}\|K_p\|_{\ca{D}}=\sqrt{1+L}\|f\|_{\ca{D}},
$$
then
$$
\|C_{\alpha}f\|_{\ca{D}}^2\leq
\left[\nu^2\left(\frac{L+2+\sqrt{L(4+L)}}{2}\right)+(1-\nu^2)(1+L)\right]\|f\|_{\ca{D}}^2,
$$
and
$\delta=\left[\nu^2\left(\frac{L+2+\sqrt{L(4+L)}}{2}\right)+(1-\nu^2)(1+L)\right]<1$
because $p\neq 0$ and $L>0$.
\end{proof}

\subsection{The essential norm}

Recall that the essential norm of an operator $T$ in a Hilbert space $\mathcal{H}$ is defined as
$\|T\|_e:=\inf\{\|T-K\|:K \mbox{ is compact }\}$, this is, the essential norm of $T$ is its norm in the
Calkin algebra. It is well known \cite{CM} that in any Hilbert space of analytic functions, we have
\begin{equation}\label{essential}
\|C_\vphi\|_e=\displaystyle\lim_n \|C_\vphi R_n\|,
\end{equation}

where $R_n$ denotes the orthogonal projection of $\mathcal{H}$ onto $z^n \mathcal{H}$.

In \cite{Sh}, it is proved that a self-map $\,\vphi:\DD\to\DD\,$ is inner if and only if the essential
norm of $C_\vphi$ in the Hardy space is equal to $\sqrt{\frac{1+|\vphi(0)|}{1-|\vphi(0)|}}.$ Because of
the analogies presented here between inner functions and univalent full-maps, one might ask: Are
full-maps characterized by the fact that the essential norm of $C_\vphi$ in the Dirichlet space is equal
to $\sqrt{\frac{L+2+\sqrt{L(4+L)}}{2}}$? where $L=\log\frac{1}{1-|\vphi(0)|^2}$. The answer is not, in
fact {\em every univalent full-map has essential norm equal to 1 in the Dirichlet space}:

\begin{theorem}
Let $\vphi:\DD\to\DD$ a univalent full-map, then $\|C_\vphi\|_e=1$ in the Dirichlet space.
\end{theorem}

\begin{proof}
Suppose first that $\vphi(0)=0$, then (\cite{MV2}) $C_\vphi$ is an isometry and equation \ref{essential}
gives:
$$\|C_\vphi\|_e=\lim_n \{\sup_{\|f\|=1} \|R_n f\|\}=\lim_n \|R_n\|=1.$$

If $\vphi(0)=p\neq 0$, then the function $\vphi_p:=\alpha_p \circ \vphi$ is a univalent full-map fixing
the origin and then for every function $f\in\mathcal{D}$ with $\|f\|=1$ we have that $\|C_\vphi R_n
f\|=\|C_{\alpha_p} R_n f\|$. Thus, $\|C_\vphi\|_e=\|C_{\alpha_p}\|_e$.

But in \cite[Cor. 5.9]{GAGR}, it is proved that the essential norm of any composition operator induced
by an automorphism of $\DD$ is equal to 1 and the result follows.
\end{proof}

\begin{center}
\textsc{Acknowledgments}
\end{center}
The authors would like to thank D. Vukoti\'{c} for suggesting the study of composition operators on the
Dirichlet space and for making available his works at his Web address.


\begin{thebibliography}{999}


\bibitem{02} E. Gallardo-Guti\'{e}rrez, and A. Montes-Rodr\'{\i}guez:
\emph{Adjoints of linear fractional composition operators on the
Dirichlet space.} Math. Ann. 327,  (2003) 117-234.
\bibitem{Bo} P. Bourdon: {\em Rudin's orthogonality problem and the
Nevanlinna counting function.} Proc. Amer. Math. Soc. 125 (1997),
1187-1192.
\bibitem{Bi} C. Bishop: {\em Orthogonal functions in $H^\infty$.}
Preprint.
\bibitem{CM} C. Cowen and B. MacCluer: {\em Composition Operators
on Spaces of Analytic Functions.} CRC Press, 1995.
\bibitem{Fe} N. Feldman: {\em Pointwise Multipliers of the Hardy space into the
Bergman space}. Illinois J. Math. 43 (1999) no. 2, 211-221
\bibitem{GAGR} G.A. Chac\'{o}n and G.R. Chac\'{o}n: {\em Some Properties of Composition Operators on the
Dirichlet Space.} Preprint. Available at the Web address:
\texttt{http://webdelprofesor.ula.ve/nucleotachira/gchacon}
\bibitem{Ha} C. Hammond: {\em The norm of a Composition Operator with Linear Fractional Symbol Acting on
the Dirichlet Space} PrePrint.
\bibitem{Ja} F. Jafari et al., editors: {\em Studies on
Composition Operators.} Comtemp. Math. Vol. 210 American Math. Soc.,
1998.
\bibitem{JM}  M. Jovovi\'{c} and B. MacCluer: {\em Composition operators on
Dirichlet spaces.} Acta Sci. Math. (Szeged) 63 (1997), 229-247.
\bibitem{MV1} M. Mart\'{\i}n and D. Vukoti\'{c}: {\em Norms and spectral radii of composition
operators acting on the Dirichlet space} To appear in  J. Math.
Anal. Appl. Available at the Web address:
\texttt{http://www.uam.es/personal\_pdi/ciencias/dragan/respub/papers.html}
\bibitem{MV2} M. Mart\'{\i}n and D. Vukoti\'{c}: {\em Isometries of the Dirichlet space among
the Composition.} To appear in Proc. Amer. Math. Soc. Available at
the Web address:
\texttt{http://www.uam.es/personal\_pdi/ciencias/dragan/respub/papers.html}
\bibitem{No} E. Nordgren: {\em
Composition operators.} Canad. J. Math. 20 (1968) 442-449.
\bibitem{Sh1} J. Shapiro: {\em Composition Operators and Classical
Function Theory.} Springer Verlag, 1993.
\bibitem{Sh} J. Shapiro: {\em What do composition operators know about inner
functions?} Monatshefte f\"{u}r Mathematik 130 (2000), 57--70.
\bibitem{Su} C. Sundberg: {\em Measures induced by analytic functions and a problem
of Walter Rudin.} J. Amer. Math. Soc. 16 (2003) 69-90.
\end{thebibliography}
\end{document}